\documentstyle[amssymb,amscd,10pt]{amsart} 
\input xy
\xyoption{all}

\newtheorem{thm}{Theorem}[section]
\newtheorem{cor}[thm]{Corollary}

\newtheorem{prop}[thm]{Proposition}
\newtheorem{lem}[thm]{Lemma}

\theoremstyle{definition}
\newtheorem{dfn}[thm]{Definition}
\newtheorem{ex}[thm]{Example}
\newtheorem{rem}[thm]{Remark}
\numberwithin{equation}{section}

\begin{document}

\newcommand{\Empieza}{$*\! >\! >$\ }
\newcommand{\Acaba}{$\ <\! <*$\ }

\newcommand{\id}{\mathop{}\mathopen{}\mathrm{id}} 
\newcommand{\Hom}{\operatorname{Hom}}
\newcommand{\Spec}{\operatorname{Spec}}
\newcommand{\C}{{\mathcal C }}
\newcommand{\Aut}{\operatorname{Aut}}
\newcommand{\Ker}{\operatorname{Ker}}
\newcommand{\Img}{\operatorname{Im}}
\newsymbol\cuadrado1003

\title[An affineness criterion for algebraic groups and applications ]
{An affineness criterion for algebraic groups and applications }
\author{C. Sancho de Salas \\ F. Sancho de Salas \\
J.B. Sancho de Salas}
\address{ \newline Carlos Sancho de Salas\newline Departamento de
Matem\'aticas\newline
Universidad de Salamanca\newline  Plaza de la Merced 1-4\\
37008 Salamanca\newline  Spain}
\email{mplu@@usal.es}
\address{ \newline Fernando Sancho de Salas\newline Departamento de
Matem\'aticas\newline
Universidad de Salamanca\newline  Plaza de la Merced 1-4\\
37008 Salamanca\newline  Spain}
\email{fsancho@@usal.es}
\address{\newline Juan B. Sancho de Salas\newline Departamento de Matem\'aticas\newline Universidad de Extremadura\newline Avda de Elvas s/n, 06006 Badajoz\newline Spain}
\email{jsancho@@unex.es}

\thanks{ The first and second authors were supported by research project MTM2017-86042-P (MINECO).
}

\begin{abstract}

We prove that a smooth and connected algebraic group $G$ is affine if and only if any invertible sheaf on any  normal $G$-variety is $G$-invariant. For the proof, a key ingredient is the following result: if $G$ is  a connected and smooth algebraic group and $\mathcal L$ is a $G$-invariant invertible  sheaf on a $G$-variety $X$, then the action of $G$ on $X$ extends  to  a projective  action on the complete linear system ${\Bbb P}(H^0(X,{\mathcal L}))$. As an application of the affineness criterion, we give a new and simple proof of the Chevalley-Barsotti Theorem on the structure of algebraic groups. 
\end{abstract}
\setcounter{section}{-1} \maketitle
\tableofcontents

\section{Introduction}

We work over an algebraically closed field. The main result of this paper is the following geometric characterization of the affineness of an algebraic group: a smooth and connected algebraic group $G$ is affine if and only if any invertible sheaf on any  normal $G$-variety is $G$-invariant. Moreover, if $G$ acts faithfully on a smooth variety $X$ and any invertible sheaf on $X$ is  $G$-invariant, then $G$ is affine (Theorem \ref{grupoafin}).
 
The first ingredient for the proof of the affineness criterion is a result (Theorem \ref{SubInv0}) that seems to be obvious but whose proof is not so simple and has its own interest: if $G$ is  a connected and smooth algebraic group and $\mathcal L$ is a $G$-invariant invertible  sheaf on a $G$-variety $X$, then the action of $G$ on $X$ extends  to  a projective  action on the complete linear system ${\Bbb P}(H^0(X,{\mathcal L}))$. Moreover, this action is coherent, that is, the complete linear system is the union of its $G$-stable  linear subsystems of finite dimension. The technical results for its proof are given in subsection \ref{fixing}, whose key point is the use of an old result of Rosenlicht about the rigidity of the group of units (Proposition \ref{Invertrosenlicht}). As a consequence of Theorem \ref{SubInv0}, we shall give (Corollary \ref{torsion}) a $G$-linearizability criterion that generalizes an essential point in  H. Sumihiro's work on $G$-equivariant completions of $G$-varieties (specifically \cite{Sumihiro2} Thm. 1.6, see also  \cite{Knop} Thm. 2.4) to (not necessarily affine) smooth and connected algebraic groups. See also Corollary \ref{Linearizability} and Remark \ref{remarklinearizable}.

A second ingredient for the proof of the affineness criterion is given by Proposition \ref{PicTrivialprop}, that states that if $Y$ is a rational normal variety, then the funtor ${\mathbf {Pic}}(Y)$, on normal schemes, is representable by a discrete and finitely generated  group scheme. This result is more or less known (see for example  Corollary 2.5 of \cite {Fossum}). We have added a proof, since it is quite elementary.

Let us comment some applications of our results. We have already mentioned the generalization to non-affine groups 
 of  Suhimiro's $G$-linearization result. A second application is an immediate proof of the finiteness of the Picard group of an affine, smooth and connected algebraic group (see \cite{Grothendieck}, p. 5-21). The third and main application is to provide a simple proof of the Chevalley-Barsotti theorem on the structure of algebraic groups (Theorem \ref{Chevalley-Barsotti}): 

{\sl If $G$ is a smooth and connected algebraic group over an algebraically closed field, then it is, in a unique way, an extension of an abelian variety   $A(G)$ by a smooth and connected affine group  $G_{\rm aff}$}. 

The first proofs of this theorem can be found in \cite{Chevalley}, \cite{Barsotti}, \cite{Rosenlicht1} or \cite{Sancho} (in the commutative case); proofs in a more modern language are given in  \cite{Conrad}, \cite{Brion}. Our proof does not require any strong result of the general theory of group schemes, except the rationality of linear groups (over an algebraically closed field); even the general theory of abelian varieties is not needed (except for the existence of the Picard variety of a smooth and complete curve); thus, it is elementary and self-contained. Leaving aside the simplicity of the proof, we obtain a new result, by characterizing the subgroup $G_{\rm aff}\subset G$ as the reduced and connected component through the origin of the subscheme of points of $G$   leaving invariant (under traslation) any invertible sheaf on $G$ (an ample invertible sheaf would suffice). The affineness  of $G_{\rm aff}$ follows from the aforementioned geometric characterization of affine groups; for the rest of the proof of the Chevalley-Barsotti theorem it only remains to prove  the representability of this group functor $G_{\rm aff}$ and the abelianity of the quotient $G/G_{\rm aff}$.

Sections 1,2  and the first part of 3, have no original results and are intended to  make a set-up of definitions and elementary results on invariant invertible sheaves that we shall need for the rest of the paper. As we have already mentioned, we work over an algebraically closed field. However, the main results are easily generalized to a perfect field by Galois descent. Beyond the perfect case, the results are false (see \cite[Prop. 5.12, Example 5.8]{Rosengarten}).

\section{Notations and definitions}\label{Notations}

We will work in the framework of \cite{EGA}. 

Throughout this paper $k$ denotes an algebraically closed field. A {\it scheme} means  a separated and noetherian scheme, locally of finite type over $k$. If, in addition, the scheme is integral, we call it a {\it variety}.  An {\it algebraic group} means a group scheme. We shall denote by ${\mathcal C}_{Sch}$ the category of schemes.

$\bullet$ Given schemes $X$ and $S$, we shall denote $X_S=S\times X$ as an $S$-scheme. For any sheaf of modules   $\mathcal M$ on $X$, we shall denote  $\mathcal M_S =\pi^*\mathcal M$, with $\pi\colon X_S\to X$ the natural projection, and for any  sheaf of modules  $\mathcal N$ on  $S$, we shall denote  ${\mathcal N}\underset k\otimes{\mathcal M}  := {\bar\pi}^*{\mathcal N} \underset{{\mathcal O}_{S\times X}}\otimes  \pi^*{\mathcal M}$, with $\bar\pi$  the projection of $S\times X$ on   $S$. In particular,  ${\mathcal M}_S={\mathcal O}_S\underset k\otimes{\mathcal M}$.

$\bullet$ For any scheme $X$, we still denote by $X$ its functor of points $X\colon {\mathcal C}_{Sch}\to \operatorname{Sets}$,   $X(S)=\operatorname{Hom}_{sch}(S,X)$. We shall use without mention Yoneda's lemma: for any functor $F\colon {\mathcal C}_{Sch}\to \operatorname{Sets}$ and any scheme $X$, giving an element of $F(X)$ is equivalent to a morphism of functors $X\to F$. If  $X=\Spec A$, we shall denote $F(A)=F(X)$. 

For any scheme $X$, $\operatorname{Pic}(X)$ denotes the set of isomorphism classes of invertible sheaves on $X$. We shall denote by ${\mathbf {Pic}}(X)$  the functor on ${\mathcal C}_{Sch}$:
$${\mathbf Pic}(X)(S)=\operatorname{Pic} (S\times X)/{\bar\pi}^*\operatorname{Pic} (S)$$ and by ${\bf Aut}_k(X)$  the functor of automorphisms of $X$: ${\bf Aut}_k(X)(S)={\rm Aut}_S(X_S)$.

A {\it pointed functor} is a pair $(F,f)$, where $F$ is a functor on ${\mathcal C}_{Sch}$ and $f$ is an element in $F(k)$. A morphism of pointed funtors $h\colon (F,f)\to (F',f')$ is a morphism of functors $h\colon F\to F'$ such that $h(k)(f)=f'$. For any pointed functor $(F,f)$ we shall denote by $F^f$ the subfunctor of $F$ defined as follows: for any scheme $S$, $F^f(S)$ is the set of morphisms $s\colon S\to F$ that factor through a morphism of schemes $S\to S'$ ($S'$ a variety) and a  morphism of pointed functors $(S',s')\to (F,f)$.  Any functor of groups $F$  is a pointed functor by taking $f=0$, the neutral element of $F(k)$, and any morphism of group functors is a  morphism of pointed functors. The corresponding subfunctor shall be denoted by $F^0$. Thus, for any scheme $X$,  ${\mathbf {Pic}}^0(X)$ is a subfunctor of ${\mathbf {Pic}}(X)$.    If   ${\mathbf {Pic}}(X) $ is representable, then  ${\mathbf {Pic}}^0(X)$  is representable by its  reduced and connected component through the origin.

$\bullet$ For any module  ${{\mathcal M}}$ on $X$ and any point $x\colon S\to X$, the  module $x^*{\mathcal M}$ on $S$  shall be called {\it fibre of ${\mathcal M}$ at $x$} and denoted by  ${\mathcal M}_x$. Thus, if $x$ is the general point of  $X$ (i.e., $x\colon X\to X$ is the identity), then ${\mathcal M}_x={\mathcal M}$. For any section  $s\in H^0(X,{\mathcal M})$, we shall denote by $s(x)$ the image of $s$ under the morphism  $ H^0(X, {\mathcal M})\to  H^0(S,x^*{\mathcal M})$ induced by $x$, and we shall say that  $s(x)$ is {\it the value of $s$ at $x$}.

$\bullet$ Let $G$ be an  algebraic group acting on a scheme $X$. The action $\tau\colon G\times X\to X$ is equivalent to giving  a morphism of group functors $\tau\colon G\to {\bf Aut}_k(X)$. For any point  $g\colon T\to G$, the corresponding automorphism  (over $T$) $X_T\to X_T$ is denoted by $\tau_g$. For any  point $x\colon T\to X$, we shall use the notation $gx=\tau_g(x)=\tau (g) (x)$. More generally, for any morphism $\tau\colon S\to {\bf Aut}_k(X)$ (i.e. an automorphism over $S$, $\tau\colon X_S\to X_S$) and any point $s\colon T\to S$, the corresponding automorphism  (over $T$) $\tau (s)\colon X_T\to X_T$ is denoted by   $\tau_s$.

\section{Rigidity of the group of units}

Let us recall a   result  due to  Rosenlicht. For any ring $A$, $A^{\rm x}$ denotes the group of units of $A$.

\begin{prop} [\cite{Rosenlicht2}, Thm. 2] \label{Invertrosenlicht} If $X$ and $S$ are varieties, then the natural morphism $$H^0(S,{\mathcal O}_{S})^{\rm {x}}\times  H^0(X,{\mathcal O}_{X})^{\rm {x}}\to
H^0(S\times X,{\mathcal O}_{S\times X})^{\rm {x}}$$ is an epimorphism with kernel $\{(\lambda,\lambda^{-1}), \lambda\in k^{\rm x}\}\simeq  k^{\rm x}$.
 
\end{prop}

\begin{cor}%[\cite{Brion2}, Lemma 2.13]\label{PicRelation} 
Let $X$ and $S$ be two varieties and $\pi\colon S\times X\to S$ the natural morphism. An  invertible sheaf ${\mathcal L}$ on  $S\times X$   is locally trivial over $S$  if and only if   $  {\mathcal L}=\pi^*\overline{\mathcal L}$ for some invertible sheaf  $\overline{\mathcal L}$ on $S$; that is,  
$${\mathbf Pic}(X)(S)=\operatorname{Pic} (S\times X)/\sim $$
where $\sim$ is the following equivalence relation: ${\mathcal L}\sim {\mathcal L}'\ \Leftrightarrow $ there exists an open covering   $\{ U_i\}$ of $S$ such that ${\mathcal L}{}_{| U_i\times X}\simeq {\mathcal L}'{}_{| U_i\times X}$.
\end{cor}

\begin{pf}  Let $\{ U_i\}$ be an open covering of  $S$ where $  {\mathcal L}$ is trivial. Let $s_0$ be a rational point of  $\cap U_i$ and let us fix an isomorphism  $\phi_{s_0}\colon  {\mathcal L}_{| \{ s_0\} \times X}\overset\sim\to {{\mathcal O}_X}$. For each $i$, let $\phi_i\colon   {\mathcal L}_{| U_i \times X}\overset\sim\to {\mathcal O}_{U_i\times X}$ be an isomorphism whose restriction to  $\{ s_0\}\times X$ is $\phi_{s_0}$. This  isomorphism is  unique up to units of $S\times X$  trivial over $\{ s_0\}\times X$, i.e., (by Proposition \ref{Invertrosenlicht}) up to units of $S$ (trivial on $s_0$). Therefore, $\phi_i,\phi_j$ differ on $U_i\cap U_j$ by a  unit   $f_{ij}$ of $U_i\cap U_j$. Obviously $\{ f_{ij}\}$ is a 1-cycle and defines an invertible sheaf  $\overline{\mathcal L}$ on $S$ such that $\pi^*\overline{\mathcal L}= {\mathcal L}$.\end{pf}

\begin{rem} This Corollary states that the Picard functor, defined as  $ {\mathbf {Pic}} (X) (S) = {\rm Pic} (X \times S) / \pi ^ * {\rm Pic} (S) $, coincides with the true Picard functor, that is, with the Zariski sheaf associated to the functor $ {Pic (X) (S)} = {\rm Pic} (X \times S) $ in the context of varieties. This justifies the given definition and allows one to use the fundamental property of $ {\mathbf {Pic}} (X) $ being  a sheaf in the Zariski topology. \end{rem}

The following result, that will be essential for the rest, follows easily  from   Proposition \ref{Invertrosenlicht}.

\begin{prop}\label{isomorHacesLineaProducto}  Let $X$ and $S$ be two  varieties and $  {\mathcal L}$ an invertible sheaf on  $S\times X$. Let  $x \in X , s \in S $ be  rational points and let us fix an isomorphism   $ {\mathcal L}_{(s,x)}=k$. If $ {\mathcal L}\simeq {\mathcal L}_1\underset k\otimes {\mathcal L}_2$ for some invertible sheaves  ${\mathcal L}_1,{\mathcal L}_2$ on $S$ and $X$, respectively, then \begin{enumerate}
\item   ${\mathcal L}_1\simeq   {\mathcal L}|_{S\times\{ x\}},\ {\mathcal L}_2\simeq    {\mathcal L}|_{\{ s\}\times X}$.
\item If  we denote $\overline{\mathcal L}_1=  {  {\mathcal L}}|_{S\times\{ x \}}, \overline{\mathcal L}_2=    {\mathcal L}|_{\{ s \}\times X}$, there exists a unique isomorphism 
$$  {\mathcal L}\overset\sim\to \overline{\mathcal L}_1\underset k\otimes \overline{\mathcal L}_2$$ which is the identity on  $\{ s \}\times X$ and $S\times\{ x \}$.
\end{enumerate}\end{prop}

\section{Invariance of invertible sheaves}

\begin{dfn}\label{defAssociated} Let $X$ be  scheme,  $\mathcal L$ an invertible sheaf on  $X$ and   $\tau\colon X_S\to X_S$ an automorphism (over $S$) (i.e., a morphism $\tau\colon S\to {\bf Aut}_k(X)$). We say that $\mathcal L$ is {\it fixed by  $\tau$}  if $\tau^*{\mathcal L}_S$ and ${\mathcal L}_S$ are isomorphic.
 
 We say that  $\mathcal L$ is {\it invariant by   $\tau $}, if there exist an invertible sheaf ${\mathcal L}^\tau$ on $S$ and an isomorphism
$$m^\tau\colon \tau^* {\mathcal L}_S\overset\sim\to  {\mathcal L}^\tau\underset k\otimes {\mathcal L}.$$ ${\mathcal L}^\tau$ is only determined up to isomorphism (Proposition \ref{isomorHacesLineaProducto}) and shall be called {\it invertible sheaf on $S$ associated to ${\mathcal L}$ and $\tau$}; the isomorphism $m^\tau$ is only determined up to units of $S\times X$ and shall be called {\it action of $\tau$ on ${\mathcal L}$}. These data will be fixed in the pointed case (subsection \ref{fixing})
\end{dfn}

\begin{rem}\label{invislocfixed} Assume that $X$  and $S$ are varieties. By Corollary \ref{PicRelation},  $\mathcal L$ is invariant by $\tau$ if and only if  $\tau^*{\mathcal L}_S$ and ${\mathcal L}_S$ are, locally on $S$ (for the Zariski topology), isomorphic; in other words, $\mathcal L$ is invariant by $\tau$ if and only if  $\mathcal L$ is, locally on $S$, fixed by $\tau$.\end{rem}

\begin{rem}\label{functorialinvariance} Sometimes, it will be useful  to enunciate the invariance of an invertible sheaf in terms of the functor of points, in the following way: an invertible sheaf $\mathcal L$ is {\it invariant by   $\tau $} if there exist an invertible sheaf $  {\mathcal L}^\tau$ on $S$ and, for any scheme $Z$ and any $Z$-valued point $(s,x)$ of $S\times X$, an isomorphism:
$${\mathcal L}_{\tau_s( x)}\overset\sim\to   ({\mathcal L}^{\tau})_{s}\!\underset{\ {\mathcal O}_Z}\otimes {\mathcal L}_x$$
which is functorial on  $(s,x)$, in the obvious sense: for any morphism $f\colon Z'\to Z$, if we denote $(s',x')=(s\circ f, x\circ f)$, then the isomorphism ${\mathcal L}_{\tau_{s'}( x')}\overset\sim\to   ({\mathcal L}^{\tau})_{s'}\!\underset{\ {\mathcal O}_{Z'}}\otimes {\mathcal L}_{x'}$ is the base change of the isomorphism ${\mathcal L}_{\tau_s( x)}\overset\sim\to   ({\mathcal L}^{\tau})_{s}\!\underset{\ {\mathcal O}_Z}\otimes {\mathcal L}_x$ by $f\colon S'\to S$.
\end{rem}

\begin{dfn} [\cite{Raynaud} V.4.1] Let $G$ be an algebraic group acting on a scheme   $X$. Let  $\tau\colon G\to {\bf Aut}_k(X)$ be the action. An invertible sheaf  $\mathcal L$ on $X$ is called {\it $G$-invariant} (respectively {\it $G$-fixed}), when it is invariant (resp. fixed) by   $\tau$. In other words (Remark \ref{functorialinvariance}),  $\mathcal L$  is $G$-invariant if and only if there exist  an invertible sheaf  $ {\mathcal L}^\tau $ on $G$ and, for any scheme $Z$ and any $Z$-valued point   $(g,x)$ of $G\times X$, a functorial isomorphism 
$$  {\mathcal L}_{gx}\overset\sim\to  ({\mathcal L}^\tau)_{g}\otimes_{{\mathcal O}_Z} {\mathcal L}_x.$$
\end{dfn}

\begin{dfn}[\cite{Mumford}, Definition 1.6]\label{G-sheaf} With the above notations, a  $G$-{\it sheaf structure} on ${\mathcal L}$ is to give, for each point $(g,x)\colon S\to G\times X$, a functorial isomophism  $\overline\tau_{(g,x)}\colon {\mathcal L}_x\overset\sim\to {\mathcal L}_{gx}$, satisfying the following associativity condition:  the triangle 
$$
\xymatrix{  {\mathcal L}_x  \ar[r]^{\overline\tau_{(g',x)} }\ar[rd]_{ \overline\tau_{(gg',x)}}& 
 {\mathcal L} _{g'x}\ar[d]^{\overline\tau_{(g,g'x)}}\\ 
  &    {\mathcal L}_{gg'x}}
$$
is commutative for any point $(g,g',x)\colon S\to G\times G\times X$.
We say that $\mathcal L$ is   $G$-{\it linearizable}  if it admits a $G$-sheaf structure.\end{dfn}

\begin{rem}\label{GsheafGfixed} The functorial isomorphism $\overline\tau_{(g,x)}\colon {\mathcal L}_x\overset\sim\to {\mathcal L}_{gx}$ is equivalent to an isomorphism ${\mathcal L}_G  \overset\sim\to \tau^*{\mathcal L} $. Hence, any $G$-linearizable invertible sheaf  is $G$-fixed.  Conversely, we shall prove   (Corollary \ref{G-lin=G-sheaf}) that if $G$ and $X$ are varieties, then any  $G$-fixed invertible sheaf on $X$ is $G$-linearizable.
\end{rem}

\begin{dfn} Let $G$ be an algebraic group and $L\colon  G\times G\to G$ (resp., $R\colon G\times G\to G$)  the action by left (resp., right) multiplication: $L(g,g')=gg'$ (resp., $R(g,g')=g'g^{-1}$). An invertible sheaf on $G$ is called {\it left-invariant} (resp., {\it right-invariant}) if it is $G$-invariant under the left action (resp., the right action) of $G$.
\end{dfn}

\begin{prop}\label{InvIzqDer} Let $G$ be an algebraic group and   ${\mathcal L}$  an invertible sheaf on  $G$. The following conditions are equivalent:

\begin{enumerate}
\item For any scheme $S$ and any $S$-valued points $g,g'$ of $G$ one has an isomorphism
$${\mathcal L}_{g  g'}\simeq {\mathcal L}_{g}\underset{{\mathcal O}_S}\otimes {\mathcal L}_{g'}$$
which is functorial on $(g, g')$.
\item  ${\mathcal L}$ is left-invariant.
\item  ${\mathcal L}$ is right-invariant.
\end{enumerate}\end{prop}

\begin{pf}  
${\mathcal L}$ is left-invariant if and only if 
there exists  an invertible sheaf  $\overline{\mathcal L} $ and a functorial isomorphism ${\mathcal L}_{g\overline g}\simeq \overline{\mathcal L}_{g}\underset{{\mathcal O}_S}\otimes {\mathcal L}_{\overline g}$. Taking $\overline g=1$, one has ${\mathcal L}_g\simeq \overline{\mathcal L}_{g}$. Thus, the left-invariance of ${\mathcal L}$ is equivalent to a functorial isomorphism:  ${\mathcal L}_{g\overline g}\simeq {\mathcal L}_{g}\underset{{\mathcal O}_S}\otimes {\mathcal L}_{\overline g}$. On the other hand, if ${\mathcal L}$ is right-invariant, there exist  an invertible sheaf  $ {\mathcal L}'$ and a functorial isomorphism ${\mathcal L}_{g\overline g^{-1}}\simeq {\mathcal L}'_{\overline g}\underset{{\mathcal O}_S}\otimes {\mathcal L}_{g}$. Taking $g=1$ one has ${\mathcal L}_{\overline g^{-1}}= {\mathcal L}'_{\overline g}$ and then   ${\mathcal L}_{g\overline g^{-1}}\simeq {\mathcal L}_{\overline g^{-1}}\otimes {\mathcal L}_{g}= {\mathcal L}_{g}\otimes {\mathcal L}_{\overline g^{-1}}$, which means that   ${\mathcal L}$ is left-invariant.\end{pf}

\begin{dfn} Since left-invariance and right-invariance of $\mathcal L$ are equivalent, we shall just say that $\mathcal L$ is {\it $G$-invariant} or {\it invariant by translation}.
\end{dfn}

\begin{prop}\label{AccionAccion0} Let ${\mathcal L}$ be an invertible sheaf on  $X$ and let $\tau\colon S\to {\bf Aut}_k(X)$, $\tau'\colon S'\to {\bf Aut}_k(X)$  be two automorphisms. Let $\tau\cdot \tau'\colon  S\times S'\to {\bf Aut}_k(X)$ be the product point, i.e., $(\tau \cdot \tau' )_{(s, s')}= \tau_s\circ  \tau'_{s'}$. If $\mathcal L$ is invariant by $\tau$ and $\tau'$, then it is invariant by  $\tau\cdot \tau'$ and one has an isomorphism 

$$ {\mathcal L}^{\tau\cdot \tau'}\overset\sim\to  {\mathcal L}^{\tau}\underset k\otimes  {\mathcal L}^{ \tau'}.$$

\end{prop}

\begin{pf}  For each $Z$-valued point  $(s,s',x)$ of $S\times S'\times X$, one has  functorial isomorphisms 
$$ {\mathcal L}_{(\tau\cdot\tau')_{(s,s')}(x)}= {\mathcal L}_{\tau_s(\tau'_{s'}(x)) } \simeq ({\mathcal L}^\tau)_s\otimes_{{\mathcal O}_Z} {\mathcal L}_{\tau'_{s'}(x)}\simeq ({\mathcal L}^\tau)_s\otimes_{{\mathcal O}_Z} ({\mathcal L}^{\tau'})_{s'}\otimes_{{\mathcal O}_Z} {\mathcal L}_x.$$
% ({\mathcal L}^{\tau \cdot\tau'})_{(s, s')}\otimes_{{\mathcal O}_Z} {\mathcal L}_x\simeq  =  &
% 
% \\ & $$ 
Then $\mathcal L$ is invariant by $\tau\cdot\tau'$ and ${\mathcal L}^{\tau \cdot\tau'}\simeq {\mathcal L}^\tau\otimes_k  {\mathcal L}^{\tau'}$.
%$({\mathcal L}^{\tau \cdot\tau'})_{(s, s')}\simeq ({\mathcal L}^\tau)_s\otimes_{{\mathcal O}_Z} ({\mathcal L}^{\tau'})_{s'}$.
\end{pf}

\begin{prop}\label{prop3.11} Let $G$ be an algebraic group acting on a scheme $X$.  If ${\mathcal L}$ is a $G$-invariant invertible sheaf on $X$, then the associated invertible sheaf  ${\mathcal L}^\tau$ on $G$ is  invariant by translation. 
\end{prop}

\begin{pf}  
By Proposition \ref{AccionAccion0}, for any $Z$-valued point $(g,g')$ of $G\times G $,  one has functorial isomorphisms $({\mathcal L}^\tau)_{gg'}= ({\mathcal L}^{\tau\cdot\tau })_{(g,g')}\simeq ({\mathcal L}^\tau)_{g }\otimes_{{\mathcal O}_Z} ({\mathcal L}^\tau)_{g' }$.
\end{pf}
\begin{rem}\label{miremark} If an invertible sheaf $\mathcal L$ on $G$ is $G$-invariant, then ${\mathcal L}^{\tau}\simeq {\mathcal L}$, because ${\mathcal L}_{g}\simeq {\mathcal L}_{g1}=({\mathcal L}^{\tau})_g\otimes {\mathcal L}_{1}\simeq ({\mathcal L}^{\tau})_g$, for each $Z$-valuated point $g\colon Z\to G$. Hence,  ${\mathcal L}$ is $G$-fixed if and only if it is trivial, because ${\mathcal O}_G\simeq {\mathcal L}^{\tau}\simeq {\mathcal L}$.
\end{rem}

\subsection{\!\!\!\! Fixing isomorphisms. The pointed case.}\label{fixing}\medskip

Let $(X,x_0)$ be a pointed scheme. A pointed invertible sheaf on $(X,x_0)$ is an invertible sheaf $\mathcal L$ on $X$ with a fixed isomorphism ${\mathcal L}_{x_0}\overset\sim\to k$.% (we shall denote  ${\mathcal L}_{x_0}= k$).

Let $\mathcal L$ be a pointed invertible sheaf on a pointed variety $(X,x_0)$. Let $(S,s_0)$ be another pointed variety and  $\tau\colon (S,s_0)\to {\bf Aut}_k(X)$ a pointed automorphism (that is, $\tau_{s_0}=\id_X$). Then $\tau^*{\mathcal L}_S$ is a pointed invertible sheaf on $(S\times X,(s_0,x_0))$, because $(\tau^*{\mathcal L}_S)_{(s_0,x_0)} \simeq {\mathcal L}_{\tau_{s_0}(x_0)} \simeq  {\mathcal L}_{x_0}\simeq k$; moreover,  $(\tau^*{\mathcal L}_S)_{\vert s_0\times X}\simeq \tau_{s_0}^* {\mathcal L}\simeq {\mathcal L}$. By Proposition \ref{isomorHacesLineaProducto}, we have:

\begin{prop}\label{pointedinvariance}   Let us denote ${\mathcal L}^\tau=(\tau^*{\mathcal L}_S)_{\vert S\times x_0}$. If ${\mathcal L}$ is invariant by $\tau$, then there exists a unique isomorphism
 \[ m^\tau\colon \tau^*{\mathcal L}_S\overset\sim\to {\mathcal L}^\tau\otimes_k {\mathcal L}\] whose restriction to
 $\{ s_0\}\times X$ and $S\times \{ x_0\}$ is the identity.
\end{prop}

\begin{rem} The invertible sheaf ${\mathcal L}^\tau$ and the morphism $m^\tau$ of Proposition \ref{pointedinvariance} are functorial in the following sense: Let $f\colon (S',s'_0)\to (S,s_0)$ be a pointed morphism, and let $\tau'=\tau\circ f$. Then $${\mathcal L}^{\tau'}=f^*{\mathcal L}^\tau$$ and the morphism $m^{\tau'}$ is the base change of $m^\tau$ by $f\times 1 \colon S'\times X\to S\times X$.
\end{rem}

Let us state Proposition \ref{pointedinvariance} in the case of the action of a group:

\begin{prop}\label{propAccion} Let  $\mathcal L$ be a pointed invertible sheaf on a pointed variety $(X,x_0)$ and  $G$ a smooth and connected algebraic group acting on $X$. Let $\tau\colon G\times X\to X$  be the action and let us denote  
$$  {\mathcal L}^\tau = (\tau^*{\mathcal L} )|_{G\times\{ x_0\}}.$$ 
 If $\mathcal L$ is $G$-invariant, then there is a unique  isomorphism  
$$m^{\tau} \colon \tau^*{\mathcal L}\overset\sim\to  {\mathcal L}^\tau \otimes_k {\mathcal L} $$
whose restriction to   $\{ 1\}\times X$ and $G\times \{ x_0\}$ is the identity. \end{prop}

\begin{dfn}\label{definaccionHazLinea}
The invertible sheaf $ {\mathcal L}^{\tau }$ of Proposition \ref{pointedinvariance} (or Proposition \ref{propAccion}) will be called  the {\it (pointed) invertible sheaf associated to $\mathcal L$ and $\tau$};  the isomorphism $m^{\tau }$   will be called the {\it (pointed) action  of $\tau$ on ${\mathcal L}$}.
\end{dfn}

\begin{rem} The pointwise description of ${\mathcal L}^\tau$ and $m^\tau$ is the following: for any $S$-valued point $g$ of $G$,  
\[ ({\mathcal L}^\tau)_g={\mathcal L}_{gx_0},\] and, for any $S$-valued point $(g,x)$ of $G\times X$, $m^\tau$ gives an isomorphism
\[   {\mathcal L}_{gx}\overset\sim\to ({\mathcal L}^\tau)_{g}\otimes_{{\mathcal O}_S} {\mathcal L}_{x}\]
which is the identity whenever $g=1$ or $x=x_0$.
\end{rem}

\begin{prop}\label{AccionAccion} Let ${\mathcal L}$ be a pointed invertible sheaf on a pointed variety $(X,x_0)$ and let $$\tau \colon (S,s_0)\to {\bf Aut}_k(X),\quad \tau' \colon (S',s'_0)\to {\bf Aut}_k(X)$$  be two pointed automorphisms   leaving   ${\mathcal L}$ invariant, with  $S,S'$ varieties. Let us denote  $$m \colon  {\mathcal L}^{\tau \cdot \tau'}\overset\sim\to  {\mathcal L}^{\tau}\underset k\otimes  {\mathcal L}^{ \tau'}$$ the unique isomorphism that restricts to 
the identity on $ \{ s_0\}\times S'$ and $S\times\{ s'_0\}$ (Proposition \ref{AccionAccion0}). The  diagram 
\begin{equation}\label{diagrama(*)}
\xymatrix{  (\tau\cdot\tau')^*{\mathcal L}_{S\times S'} = \ar[d]_{m^{\tau\cdot\tau'}}&\hskip -1cm
(\tau')^*[(\tau^*{\mathcal L}_{S})_{ S'}] \ar[r]^{\hskip -.5cm m^{\tau  }}& 
 {\mathcal L}^{\tau  }\underset {k}\otimes  (\tau')^*{\mathcal L}_{S'} \ar[d]^{1\otimes m^{ \tau' }}\\ 
 {\mathcal L}^{\tau\cdot\tau' }\underset {k}\otimes {\mathcal L}
  \ar[rr]^{m \otimes 1\ \ \ \ } & &   {\mathcal L}^{\tau }\underset {k}\otimes {\mathcal L}^{\tau' }\underset {k}\otimes {\mathcal L}}
\end{equation}
is commutative.
\end{prop}

\begin{pf} Both compositions are   isomorphisms   $ (\tau \cdot\tau')^*{\mathcal L}_{S \times S'}\to  {\mathcal L}^{\tau  }\underset {k}\otimes {\mathcal L}^{\tau' }\underset {k}\otimes {\mathcal L}$ that restrict to the identity on  $\{ s_0 \}\times \{ s_0'\}\times X , \{ s_0\}\times S'\times \{ x_0\}$ and $S\times \{  s_0'\}\times \{ x_0\}$, so they coincide by Proposition \ref{isomorHacesLineaProducto}.
\end{pf}

Taking the fiber at a $Z$-valued point $(s,s',x)$ of $S\times S'\times X$ in \ref{diagrama(*)} we obtain the commutative diagram:
\begin{equation}\label{diagrama(**)}
\xymatrix{  {\mathcal L}_{\tau_s(\tau'_{s'}(x))}  \ar[r]^{m^{\tau_s }} \ar[d]_{ m^{ (\tau\cdot\tau')_{(s,s')} }}& 
( {\mathcal L}^{\tau_s  })_s\underset {{\mathcal O}_Z}\otimes  {\mathcal L} _{\tau'_{s'}(x)}\ar[d]^{1\otimes m^{ \tau' }}\\ 
 ({\mathcal L}^{\tau\cdot\tau' })_{(s,s')}\underset {{\mathcal O}_Z}\otimes {\mathcal L}_x
  \ar[r]^{m _{(s,s')}\otimes 1\ \ }  &   ({\mathcal L}^{\tau })_s\underset {{\mathcal O}_Z}\otimes ({\mathcal L}^{\tau' })_{s'}\underset {{\mathcal O}_Z}\otimes {\mathcal L}_x}
\end{equation}

\begin{cor}[\cite{Knop} Lemma 2.3, \cite{Brion2} Lemma 2.9]\label{G-lin=G-sheaf}  Let $G$ be a  smooth and connected algebraic group acting on a variety $X$ and let  $\mathcal L$ be an invertible sheaf on $X$.  Then $\mathcal L$ is $G$-linearizable if and only if it is  $G$-fixed.\end{cor}

\begin{pf} If $\mathcal L$ is $G$-linearizable, then it is $G$ fixed (Remark \ref{GsheafGfixed}). Conversely, let  $\mathcal L$ be a $G$-fixed invertible sheaf on $X$. Let us fix a rational point $x_0\in X$ and an isomorphism ${\mathcal L}_{x_0}=k$. Since $\mathcal L$ is $G$-fixed, it is $G$-invariant and there is an isomorphism $\phi\colon {\mathcal L}^\tau \simeq{\mathcal O}_G$.  Substituting in diagram \ref{diagrama(**)} we obtain 

$$\xymatrix{  {\mathcal L}_{gg'x}  \ar[r]^{m^{\tau _g}} \ar[d]_{ m^{ \tau_{gg'} }}& 
 {\mathcal O}_g\underset {{\mathcal O}_Z}\otimes  {\mathcal L} _{g'x}\ar[d]^{1\otimes m^{ \tau_{g'} }}\\ 
 {\mathcal O}_{gg'}\underset {{\mathcal O}_Z}\otimes {\mathcal L}_x
  \ar[r]^{m \otimes 1\ \ }  &   {\mathcal O }_g\underset {{\mathcal O}_Z}\otimes {\mathcal O }_{g'}\underset {{\mathcal O}_Z}\otimes {\mathcal L}_x}.
$$
Taking $g'=1$, one obtains a  functorial isomorphism   $\overline\tau_{(g,x)}=(m^{\tau_g})^{-1}\colon  {\mathcal L }_x\overset\sim\to  {\mathcal L }_{gx}$, i.e.  a  $G$-sheaf structure on ${\mathcal L}$ (Definition \ref{G-sheaf}).\end{pf}

\section{Projective actions on sections of invariant invertible sheaves}{}

Let $\mathcal E$ be a quasi-coherent module on a scheme  $S$. We shall denote  $${\Bbb P}_S({\mathcal E})={\bf Proj}_S\, S^{\cdot }_{{\mathcal O}_S}{\mathcal E}$$ 
When $S=\Spec k$ and ${\mathcal E}=E$ is a $k$-vector  space,  we shall just denote it by  ${\Bbb P}( E)$. The scheme ${\Bbb P}_S({\mathcal E}) $ represents the functor on ${\mathcal C}_{\rm Sch}$ that assigns to each scheme $T$ the  pairs $(s, \omega )$ where  $s\colon T\to S$ is a point of $S$ and  $\omega\colon  {\mathcal E}_s\to \widetilde {\mathcal L}$ is an invertible quotient of  ${\mathcal E}_s$. 

An isomorphism $T\colon {\mathcal E}\to {\mathcal E}'$ induces an isomorphism $[T]\colon {\Bbb P}_S({\mathcal E}')\to {\Bbb P}_S({\mathcal E})$ of $S$-schemes. For any base change $f\colon S'\to S$ one has a natural isomorphism ${\Bbb P}_{S'}({f^*\mathcal E})={\Bbb P}_S({\mathcal E})\times_SS'$. Finally, we shall denote ${\Bbb P}^*_S({\mathcal E})={\Bbb P}_S({\mathcal E}^*)$, where  ${\mathcal E}^*={\mathcal H}om_{{\mathcal O}_S}({\mathcal E},{\mathcal O}_S)$ is the dual sheaf.

\begin{dfn}\label{projectivity} Let $\mathcal E$ be a quasi-coherent module on   $S$. An automorphism of $S$-schemes  
$$\tau\colon {\mathbb P}_S({\mathcal E})\to {\mathbb P}_S({\mathcal E})$$ is called a {\it projectivity} if, 
Zariski locally on  $S$, it is the automorphism induced by an automorphism of   $\mathcal O_S$-modules, $T\colon {\mathcal E}\to {\mathcal E}$.  
\end{dfn}

\begin{ex}\label{projectivityexample} Let $S$ be a scheme,   ${\mathcal E}$  a quasi-coherent module on $S$ and ${\mathcal L}$   an invertible sheaf on $S$. One has a canonical isomorphism ${\mathbb P}_S({\mathcal E})= {\mathbb P}_S({\mathcal E} \otimes_{{\mathcal O}_S} {\mathcal L})$. Hence, an isomorphism  $m\colon {\mathcal E}\to {\mathcal L} \otimes_{{\mathcal O}_S} {\mathcal E}$ induces an automorphism of $S$-schemes:
$$ [m]\colon {\mathbb P}_S({\mathcal E})\to {\mathbb P}_S({\mathcal E})$$
which  is a projectivity. 
\end{ex}

\begin{dfn} \label{definicionCoherente}  
(1) Let   $T\colon {\mathcal E}\to {\mathcal E}$ be an endomorphism of modules. It induces an ${\mathcal O}_S[X]$-module structure   on $\mathcal E$ by $X\cdot e=T(e)$. We say that $T$ is {\it coherent} if, locally on $S$, 
${\mathcal E}$  is, as an   ${\mathcal O}_S[X]$-module, the direct limit of its sub-${\mathcal O}_S[X]$-modules that are coherent as  ${\mathcal O}_S$-modules. That is, for each affine open subset $\Spec A=U$ of $S$ and each section  $s\in {\mathcal E}(U)$, $A[T]\cdot s$ is  a finitely generated $A$-module. 

(2) A projectivity $\tau\colon {\mathbb P}_S({\mathcal E})\to {\mathbb P}_S({\mathcal E})$ is called {\it coherent}  if, locally on $S$, it is the projectivity associated to a coherent automorphism of ${\mathcal E}$ .
  
\end{dfn}

For the rest of this section,  $X$ is a  scheme and ${\mathcal L}$ is an invertible sheaf on  $X$. Let us denote   $E =H^0(X,{\mathcal L})$. For each scheme $S$, let $\pi\colon S\times X\to S$ be the projection. One has a natural isomorphism $E_S=\pi_*({\mathcal L}_S)$. Let ${\mathcal B} \subset X$ be the locus of base points of  $E $, i.e., the locus of points of $X$ where the morphism  $\phi\colon {\mathcal O}_X\underset k\otimes E \to{\mathcal L}$ is not surjective. It is a closed subscheme whose ideal of definition is the $0$-th Fitting ideal of the cokernel of $\phi$. The morphism $\phi$ induces a natural morphism
$$\nu \colon X-{\mathcal B} \to {\Bbb P} (E )$$
that maps a rational point  $x\in X-{\mathcal B}$ into the quotient 
 $$\nu_x\colon E\to {\mathcal L}_x$$
  defined by $\nu_x(s)=s(x)$,  the value of $s$ at $x$.
 If $x\colon S\to X-{\mathcal B}$ is an arbitrary point, one obtains analogously the quotient   $\nu_x\colon E_S \to {\mathcal L}_x$.

\begin{prop}\label{accionSecciones} Let $X$ be a scheme,   ${\mathcal L}$ an invertible sheaf on $X$ and $E=H^0(X,{\mathcal L})$.   If $\tau\colon X_S\to X_S$ is an automophism (over $S$) that leaves   ${\mathcal L}$ invariant, then  there exists a projectivity  $\overline \tau\colon {\mathbb P} (E)_S\to {\mathbb P} (E)_S $ which is compatible with $\tau$, i.e., the diagram 
$$\xymatrix{ (X -{\mathcal B})_S \ar[r]^{\tau}\ar[d]_{\nu_S} & (X-{\mathcal B})_S\ar[d]^{\nu_S}\\
{\mathbb P} (E)_S \ar[r]^{\overline\tau } & {\mathbb P} (E)_S
}$$ 
is commutative. 
\end{prop}

\begin{pf}  Since $\mathcal L$ is invariant by $\tau$, there exists an isomorphism  $m^\tau\colon \tau^*{\mathcal L}_S\overset\sim\to {\mathcal L}^\tau\underset k\otimes {\mathcal L}$. Taking the direct image  by  $\pi\colon S\times X\to S$, one has an isomorphism  
$$\Gamma^\tau\colon E_S=\pi_*({\mathcal L}_S) =\pi_*(\tau^*{\mathcal L}_S)\overset \sim\longrightarrow  {\mathcal L}^\tau \otimes_k E  =  {\mathcal L}^\tau \otimes_{{\mathcal O}_S} E_S $$
hence a projectivity (see Example \ref{projectivityexample}):
$$\overline\tau =[\Gamma^\tau]\colon{\mathbb P} ({\mathcal L}^\tau \otimes_{{\mathcal O}_S} E_S)={\mathbb P} (E)_S \overset\sim\longrightarrow {\mathbb P} (E)_S.$$  

For each $Z$-valuated point $(s,x)\colon Z\to S\times X$, one has the commutative diagram (we shall denote by  $\pi\colon X_Z\to Z$ the natural projection) :
$$\xymatrix {   E_Z  \ar[d]_{\nu_{ \tau_s(x)}} & \hskip -40pt = \\
 {\mathcal L}_{ \tau_s(x)} &\hskip -40pt  =
} 
\hskip -20pt\xymatrix{  \pi_*(\tau_{ s}^*{\mathcal L}_{ Z})\ar[d]_{\nu_{x}} 
\ar[r]^{ (\Gamma^\tau)_s \ \ \ }&   ({\mathcal L}^\tau)_{ s} \otimes_{{\mathcal O}_Z} E_Z \ar[d]^{1\otimes \nu_{x}}\\
 (\tau_{ s}^*{\mathcal L}_{ Z})_{x}
\ar[r]^{
(m^\tau)_{(s,x)
} \ \ \ } &   ({\mathcal L}^\tau)_{ s}\otimes_{{\mathcal O}_Z}{\mathcal L}_{x}
}$$
$$\xymatrix{  & &}$$
hence the diagram 
$$\xymatrix{   (X -{\mathcal B})_S\ar[r]^{\tau } \ar[d]_{\nu_S }& (X  -{\mathcal B})_S  \ar[d]^{\nu_S }\\ {\Bbb P} (E)_{ S} \ar[r]^{\overline\tau } & {\Bbb P}(E)_{ S}}$$
is commutative.  
\end{pf}

The projectiviy $\overline \tau$ is not unique in general. However, one has

\begin{prop}\label{accionSecciones2} Let $X$ be a variety,   ${\mathcal L}$ an invertible sheaf on $X$ and $E=H^0(X,{\mathcal L})$. Let $(S,s_0)$ be a pointed variety and $\tau\colon (S,s_0)\to {\bf Aut }(X)$ a pointed morphism (i.e.,  $\tau_{s_0}=\id_X$).
If ${\mathcal L}$ is invariant by $\tau$, there exists a unique projectivity $$\overline \tau\colon {\mathbb P} (E)_S\to {\mathbb P} (E)_S $$ which is compatible with $\tau$ (as in Proposition \ref{accionSecciones}) and such that   $\overline\tau_{s_0}=\id_{{\mathbb P} (E)}$. Moreover, the assignement  $\tau\mapsto \overline \tau$ satisfies:
\begin{enumerate}

\item It is functorial on $(S,s_0)$: for any pointed morphism $ (S',s'_0)\to (S,s_0)$, one has $\overline{\tau_{S'}} = {\overline\tau}_{S'}$, where $\tau_{S'}$ (resp. ${\overline\tau}_{S'}$) is the base change of $\tau$ (resp. $\overline\tau$) by $S'\to S$.

\item It is compatible with compositions:  if $\tau_1 , \tau_2\colon (S,s)\to {\bf Aut }(X)$ are two pointed morphisms leaving  ${\mathcal L}$ invariant, then 
$$\overline{\tau_1} \circ\overline{\tau_2} = \overline{\tau_1\circ\tau_2 }.$$ 
\end{enumerate}
\end{prop}

 \begin{pf} Let us fix a rational point $x_0\in X$ and an isomorphism ${\mathcal L}_{x_0}=k$. Let ${\mathcal L}^\tau$ and $m^\tau\colon \tau^*{\mathcal L}_S\to {\mathcal L}^\tau\otimes_k {\mathcal L}$ be those of Proposition \ref{pointedinvariance}. Let $\Gamma^\tau\colon E_S\to {\mathcal L}^\tau\otimes_{{\mathcal O}_S} E_S$ be the induced isomorphism (as in the proof of Proposition \ref{accionSecciones}) and $\overline\tau=[\Gamma^\tau]$ the associated projectivity. We know that $\overline\tau$ is compatible with $\tau$; moreover,  $\overline\tau_{s_0}$ is the identity because $m^\tau$ is the identity on $s_0\times X$.

Let us prove that $\overline\tau$ is unique. It suffices to prove that if $\overline \tau$ is a projectivity of   ${\mathbb P} (E)_S$ compatible with $\id_X$ (Proposition \ref{accionSecciones}) and  $\overline\tau_{s_0}=\id_X$, then $\overline \tau$ is the identity projectivity. Let us denote $\overline X=X-\mathcal B$. By the compatibility condition we have $\overline \tau_{| \nu_S (\overline X_S)}=\id_{\nu_S (\overline X_S)}$. Now, we may assume that $S$ is an affine variety $S=\Spec A$ and that $\overline \tau$  is induced by an isomorphism  $L\colon E_S\overset\sim\to E_S$.  Let   $\theta\in A\otimes_k  E$ be part of a basis. In $\overline X_S$ one has the schematic equality  $(\theta)_0=\nu_S^{-1}(H_\theta\,\cap \nu_S(\overline X_S))$, where $H_\theta$ is the hyperplane defined by  $\theta$. Since   $H_{L(\theta)}=\overline \tau (H_\theta)$, it follows that $(L(\theta))_0=(\theta)_0$. This means that  $L(\theta)$ and $\theta$ differ by a unit $a\otimes f$ of $\overline X_S=S\times \overline X$. Since $\overline\tau_{s_0}=\id$, we may assume that $f=1$, thus $<\theta>=<L(\theta)>$ as sub-$A$-modules of $A\otimes_kE$; hence $\overline\tau=\id$.

The functoriality of  $\tau\mapsto \overline \tau$ and its compatibility with compositions follow immediately from the uniqueness of $\overline\tau$. 
 \end{pf}

\begin{thm}{\label{SubInv0}} Let $G$ be a smooth and connected algebraic group acting on a variety  $X$ and let  ${\mathcal L}$ be a $G$-invariant invertible sheaf on  $X$. The action of  $G$ on $X$ extends to an action of $G$ on  ${\mathbb P} (H^0(X,{\mathcal L}))$ by coherent projectivities.  \end{thm}

\begin{pf} We shall denote $E=H^0(X,{\mathcal L})$. By  Proposition \ref{accionSecciones2}, the action $\tau\colon X_G\to X_G$ extends to a projectivity $\overline\tau\colon {\mathbb P}(E)_G\to {\mathbb P}(E)_G$ (unique such that ${\overline\tau}_1=\id$).
Thus we have a (pointed) morphism $G\to {\bf Aut} ({\mathbb P}(E))$. It remains to prove that this is a morphism of group functors and that $\overline\tau$ is coherent.

(1) $G\to {\bf Aut} ({\mathbb P}(E))$ is a morphism of group functors. This is a consequence of (1) and (2) of Proposition \ref{accionSecciones2}. Indeed, let $\tau_1, \tau_2\colon X_{G\times G}\to X_{G\times G}$ be defined by $\tau_1(g,g',x)=(g,g',gx)$,  $\tau_2(g,g',x)=(g,g',g'x)$. By (2) of  Proposition \ref{accionSecciones2}, one has $\overline{\tau_1\circ\tau_2} = \overline{\tau_1}\circ\overline{\tau_2}$. Moreover, $\tau_1$ (resp. $\tau_2$, $\tau_1\circ\tau_2$) is the base change of $\tau$ by the first projection $G\times G\to G$ (resp. the second projection, the multiplication morphism). Hence, for any scheme $T$ and any $T$-valued point $(g,g')$ of $G\times G$, one has
\[ {\overline\tau}_{gg'}=({\overline{\tau_1\circ\tau_2}})_{(g,g')}= {\overline{\tau_1}}_{(g,g')} \circ {\overline{\tau_2}}_{(g,g')} = {\overline\tau}_g\circ {\overline\tau}_{g'}.\]

(2) $\overline\tau$ is coherent. Let us fix a rational point $x_0\in X$ and an isomorphism ${\mathcal L}_{x_0}=k$. Let ${\mathcal L}^\tau$ and $m^\tau\colon \tau^*{\mathcal L}\to {\mathcal L}^\tau\otimes_k\mathcal L$ be those of Proposition \ref{propAccion}. Then one has an isomorphism $\Gamma^\tau\colon E_G\to {\mathcal L}^\tau\otimes_k E$ and $\overline\tau =[\Gamma^\tau]$. If we denote $V=H^0(G,{\mathcal L}^\tau)$, $\Gamma^\tau$ is equivalent to a morphism (that we still denote $\Gamma^\tau$)
\[\Gamma^\tau\colon E\to V\otimes_k E\]

Let $e\in E$. We have to prove that there exists a finite vector subspace $E'\subseteq E$ containing $e$ and such that $\Gamma^\tau(E')\subseteq V\otimes_k E'$. 

Let $\{ v_i\}$ be a basis of $V$ and let $e_i\in E$ be such that
\[ \Gamma^\tau  (  e) =\sum_i v_i\otimes e_i.\]
Let $E'$ be the subspace generated 
by the vectors $ e_i$, which is finite dimensional. 

Let us see that $e\in E'$. Let $V\to ({\mathcal L}^\tau)_1=k$ be the evaluation at $1$, $v\mapsto v(1)$. The composition 
\[ \aligned E\overset{\Gamma^\tau}\to V\otimes_k E &\to E\\ v\otimes e'&\mapsto v(1)e'\endaligned\] is the identity, because $\tau_1=\id$. Hence,  $e=\sum_i v_i(1)e_i$, so $e\in E'$.

It remains to prove that $E'$ is stable by $\Gamma^\tau$. This is a consequence of Proposition \ref{AccionAccion}. Let us give the details:

Let us consider the automorphisms $\tau_1=\tau_2=\tau$ and the product $\tau_1\cdot\tau_2\colon X_{G\times G}\to X_{G\times G}$. Let us denote 
\[\aligned 1\otimes\Gamma^\tau\colon V\otimes_kE&\to V\otimes_kV\otimes_kE  \\ v\otimes e' &\mapsto v\otimes \Gamma^\tau( e') 
\endaligned\]
On the other hand, let us consider the isomorphism $ m^{\tau_1\cdot\tau_2}\colon (\tau_1\cdot\tau_2)^*{\mathcal L}_{G\times G}\to {\mathcal L}^{\tau_1\cdot\tau_2} \otimes_k {\mathcal L}$ and let $$\Gamma^{\tau_1\cdot\tau_2}\colon  E\to W\otimes_k E$$ be the induced  morphism, where $W=H^0(G\times G,{\mathcal L}^{\tau_1\cdot\tau_2})$. The isomorphism $m\colon {\mathcal L}^{\tau_1\cdot\tau_2}\to {\mathcal L}^\tau\otimes_k  {\mathcal L}^\tau$ of Proposition \ref{AccionAccion} induces an isomorphism $m \colon W\to V\otimes_k V$, hence an isomorphism $m \otimes 1\colon W\otimes_k E\to V\otimes_k V\otimes_k E$. By  Proposition \ref{AccionAccion}, one has  
\[ (1\otimes\Gamma^{\tau})\circ  \Gamma^{\tau} = (m \otimes 1)\circ \Gamma^{\tau_1\cdot\tau_2}.\]
Now, 
\begin{equation}\label{eq1} ((1\otimes\Gamma^{\tau})\circ  \Gamma^{\tau})( e)=\sum_i v_i\otimes \Gamma^\tau( e_i).\end{equation} 
On the other hand,  $\tau_1\cdot\tau_2$ is the automorphism obtained from $\tau$ under the base change $\phi\colon G\times G\to G$, $\phi(g,g')=gg'$; hence, ${\mathcal L}^{\tau_1\cdot\tau_2}=\phi^*{\mathcal L}^\tau$, which induces a morphism $\phi^*\colon V\to W$, and  $\Gamma^{\tau_1\cdot\tau_2}=(\phi^*\otimes 1)\circ \Gamma^\tau$. Hence 
\[\Gamma^{\tau_1\cdot\tau_2} (  e)=\sum_t w_t\otimes e_t,\quad w_t=\phi^*(v_t).\]   If we write $ m   (w_t)=\sum_{j,k} \lambda^t_{jk} v_j\otimes v_k$, then 
\begin{equation}\label{eq2} ((m\otimes 1)\circ \Gamma^{\tau_1\cdot\tau_2}) ( e)=\sum_{j,k,t}\lambda^t_{jk} v_j\otimes v_k \otimes e_t
\end{equation}

From equations \eqref{eq1} and \eqref{eq2}, we obtain that $\Gamma^\tau ( e_i)= \sum_{k,t}\lambda_{ik}^t  v_k\otimes e_t$ and then $\Gamma^\tau (  E')\subseteq V\otimes_k E'$.  
\end{pf}

\begin{cor} \label{torsion} Let ${\mathcal L}$  be a  $G$-invariant invertible sheaf on a $G$-variety $X$. If  the  sheaf ${\mathcal L}^{\tau}$ on $G$ associated to ${\mathcal L}$ (Definition \ref{defAssociated}) is effective then there exists $n>0$ such that ${\mathcal L}^n$ is $G$-linearizable.
\end{cor}

\begin{pf} (a) Let us prove the statement when ${\mathcal L}$ is effective (instead of ${\mathcal L}^{\tau}$). 
 Let us take $s_1,\dots , s_r\in H^0(X,{\mathcal L})$ such that  $(s_1)_0\cap\cdots\cap (s_r)_0=\mathcal B$. Let $V\subset H^0(X,{\mathcal L})$ be a finite dimensional subspace containing $\{ s_1,\dots ,s_r\}$ and such that ${\Bbb P}(V)$ is stable under the action of $G$ on ${\Bbb P}(H^0(X,{\mathcal L}))$   (Theorem \ref{SubInv0}). We have a $G$-morphism $\pi\colon X-{\mathcal B}\to {\Bbb P}(V)$. Let  $n=\operatorname{dim}_kV$. Since ${\mathcal O}_{{\Bbb P}(V)}(n)\simeq \Lambda^{n-1}_{{\mathcal O}_{{\Bbb P}(V)}}\Omega_{{\Bbb P}(V)/k}^*$, we obtain that ${\mathcal O}_{{\Bbb P}(V)}(n)$ is   ${\mathbf PGl}_k(V)$-linearizable, so it is   $G$-linearizable and then ${\mathcal L}^n{}_{|X-{\mathcal B}}=\pi^*{\mathcal O}_{{\Bbb P}(V)}(n)$ is   $G$-linearizable. Since ${\mathcal L}^{\tau}=({\mathcal L}_{|X-\mathcal B})^{\tau}$, one has  $({\mathcal L}^n)^{\tau}=({\mathcal L}_{|X-\mathcal B}^n)^{\tau}\simeq {\mathcal O}_G$; thus ${\mathcal L}^n$ is fixed by $G$ and we conclude by Corollary \ref{G-lin=G-sheaf}.
 
Now, let us prove the statement. Since ${\mathcal L}^{\tau}$ is   effective and $G$-invariant (Proposition \ref{prop3.11}),   $({\mathcal L}^{\tau })^n$  is $G$-linearizable by (a). Then,   $({\mathcal L}^{\tau })^n$ is  trivial (Remark \ref{miremark}), so ${\mathcal L}^n$ is $G$-fixed. By Corollary \ref{G-lin=G-sheaf}, ${\mathcal L}^n$ is $G$-linearizable.\end{pf}

\begin{rem} 
As has been shown in the proof, Corollary \ref{torsion} still holds if we assume that the $G$-invariant sheaf $\mathcal L$ is   effective (instead of ${\mathcal L}^\tau$). But this is a stronger condition, since the effectiveness of $\mathcal L$ implies that of ${\mathcal L}^\tau$; indeed, taking global sections in the isomorphism $\tau^*{\mathcal L}\overset\sim\to {\mathcal L}^{\tau}\underset k\otimes {\mathcal L}$ one has $$0\neq H^0(G,{\mathcal O}_G)\underset k\otimes H^0( X,{\mathcal L})=H^0(G\times X,\tau^*{\mathcal L})\simeq H^0(G,{\mathcal L}^{\tau})\underset k\otimes H^0( X,{\mathcal L})$$ and then $H^0(G,{\mathcal L}^{\tau})\neq 0$. \end{rem}

\begin{rem} The effectiveness of ${\mathcal L}^{\tau}$ cannot be eliminated in Corollary \ref{torsion}. In fact, let $G$ be an abelian variety and $G^*={\mathbf {Pic}}^0(X)$ its dual abelian variety. The invariant sheaves on $G$ are the rational points of $ G^* $. These  sheaves satisfy that ${\mathcal L}^{\tau }={\mathcal L}$ (Remark \ref{miremark}). Moreover ${\mathcal L}$ is linearizable if and only if it is trivial. Therefore ${\mathcal L}^n$ is linearizable for some $n>0$ if and only if ${\mathcal L}$ is a torsion point of $G^*$.  \end{rem}

\section{Representability of the Pic funtor on rational varieties}

We are going to prove a functorial version of a well-known elementary result (see Corollary 2.5 of \cite {Fossum}). This is:

\begin{prop} \label{PicTrivialprop} If $X$ is a normal and rational variety, then the funtor ${\mathbf {Pic}}(X)$, on the category of  normal schemes,  is representable by the constant group scheme associated to a finitely-generated abelian group. This holds in particular  if $X$ is an affine, smooth and connected algebraic group.
\end{prop}

\begin{pf} We shall need some preliminary considerations:

Let $ X $ be a normal variety and $ Y \subset X $ a closed subset. Let $ D_1, \dots, D_r \subset Y $ be the codimension 1 (in $X$) irreducible components of $ Y $. Let us denote by $ Div (X) $ the group of Weil divisors of $X$, $ Div_Y (X) $  the free subgroup generated by $ D_1, \dots, D_r $ and $ Car_Y (X) \subset Div_Y (X) $ the finitely generated subgroup  of Cartier divisors.  Now let $ {\rm Pic} _Y (X) $ be the quotient of $ Car_Y (X) $ by  linear equivalence. In other words:
$$\aligned {\rm Pic}_Y(X)=\{&\text{Group} \text{ of  (isomorphism classes of) invertible sheaves on $ X $ }\\
&\text{that are associated  to some Cartier divisor  supported on $ Y $}\}\endaligned$$
We denote by ${\mathbf Pic}_Y(X)$ the constant  sheaf associated to ${\rm Pic}_Y(X)$, that is:
 $${\mathbf {Pic}}_Y(X)[S]= {{\rm Pic}}_Y(X)$$ for each variety $S$. 

Obviously  $ {\mathbf {Pic}}_Y(X)$ is representable by  the scheme whose underlying space is the finitely generated commutative group  $ {\rm Pic} _Y (X) $ (with the discrete topology).

If we denote $U=X-Y$,  one has the  exact sequence:
$$0\to {\rm Pic}_Y(X)\to {\rm Pic}(X)\overset {i^*}\to  {\rm Pic}(U)$$
where $ i ^ * $ is the restriction morphism.

\begin{lem}\label{sucExac}

With the above notation, one has an exact sequence of group funtors (over normal - non necessarily connected - schemes):
$$0\longrightarrow {\mathbf Pic}_Y(X)\longrightarrow {\mathbf Pic}(X)\overset {i^*}\longrightarrow  {\mathbf Pic}(U).$$
Moreover, if $ X $ is smooth, then $ i ^ * $ is surjective over  the category of smooth schemes.
\end{lem}

\begin{pf} For any normal variety $S$ one has an exact sequence   $$0\to {\rm Pic}_{Y_S}(X_S)\to {\rm Pic}(X_S)\overset {i^*_S}\to  {\rm Pic}(U_S)$$ and 
since $\pi_X^*Car_Y (X)=Car_{Y_S} (X_S)$, where $\pi_X\colon X\times S\to X$ is the natural projection, one has an isomorphism  ${\rm Pic}_{Y}(X)\overset{\pi_X^*}\simeq {\rm Pic}_{Y_S}(X_S)$. The conclusion  follows.

The surjectivity of $i^*$ in the smooth case  is a consequence of the fact that  the  restriction morphism  on 
Weil divisors is surjective and that every 
Weil divisor is a Cartier divisor.
\end{pf}

Now let us prove  Proposition \ref{PicTrivialprop}:

To see the representability of $ {\mathbf Pic} (X) $ by a discrete scheme it   suffices to prove  that ${\rm Pic}(X\times S)={\rm Pic}(X)\times {\rm Pic}(S)$  for each  normal variety $ S $.

Let us denote by $ X_ {sing} $ the closed subset of singular points of $X$ and $ X_ {reg} = X-X_ {sing} $. Since $ X \times S $ is normal, $ (X \times S) _ {sing} $ does not contain Weil divisors and therefore the restriction morphism $ {\rm Pic} (X \times S) \to {\rm Pic} (X_ {reg} \times S_ {reg}) $ is injective. Hence, an invertible sheaf on $ X \times S $ is isomorphic to $ {\mathcal L} _1 \underset k \otimes {\mathcal L} _2 $ for some invertible sheaves $ {\mathcal L} _1 , {\mathcal L} _2 $ on $X$ and $S$, if and only if it is so on $ X_ {reg} \times S_ {reg} $ (see Proposition \ref{isomorHacesLineaProducto}). Thus, replacing $ X $ and $ S $ by $ X_ {reg} $ and $ S_ {reg} $, we can assume that $ X $ and $ S $ are smooth.

Since $X$ is rational, there exists an open subset $ U$ of $X$ which is isomorphic to an open subset $ V $ of an affine space $ {\mathbb A}^n $.
It is well known (see \cite{Hartshorne}, Chap. II, Proposition 6.6) that $ {\rm Pic} ({\mathbb A}^n \times S) = {\rm Pic} (S) $ for each smooth variety $ S $ and then $ {\mathbf Pic} ({\mathbb A}^n) = 0 $. Since $ {\mathbb A}^n $ is smooth, we obtain that $ {\mathbf Pic} (V)= 0 $  by  Lemma \ref{sucExac}. Again, applying Lemma \ref{sucExac} to $X$ and $U$,  where $Y=X- U$, we obtain that   $ { \mathbf Pic} (X) = {\mathbf Pic} _Y (X) $, because   $ {\mathbf Pic} (U) \simeq {\mathbf Pic} (V) = 0 $. We conclude because ${\mathbf Pic} _Y (X) $  is representable by a discrete and finitely generated group scheme.\end{pf}

\medskip\section{The affineness criterion for algebraic groups.}\medskip

The aim of this section is to prove the following criterion for affineness of a group:
 
\begin{thm} \label{grupoafin} Let $G$ be a smooth and connected algebraic group.\begin{enumerate}

\item If  $G$ is affine, then any invertible sheaf on any normal $G$-variety is  $G$-invariant. 

\item Conversely, if $G$ acts faithfully on a   smooth variety $X$ and for each in\-ver\-ti\-ble sheaf  ${\mathcal L}$  on $X$, there exists some $n> 0$ such that ${\mathcal L}^n$ is $G$-invariant, then $G$ is affine. 

\end{enumerate}\end{thm}

\begin{pf} (1) Let ${\mathcal L}$ be an invertible sheaf on a normal $G$-variety  $X$ and let us denote  by  $\tau\colon G\times X\to X$, the action. By Yoneda's Lemma, the invertible sheaf $\tau^*{\mathcal L}$ on $G\times X$ defines a morphism $\varphi_{\mathcal L}\colon X\to {\mathbf {Pic}}(G)$ that maps a point $x\colon S\to X$ to $  (1\times x)^*(\tau^*{\mathcal L})\in{\rm Pic}(G\times S)$. Since $G$ is rational, ${\mathbf {Pic}}(G)$ is a discrete scheme (Proposition \ref{PicTrivialprop}) and then $\varphi_{\mathcal L}$ is constant, because $X$ is connected. Thus, $\tau^*{\mathcal L}\simeq {\mathcal L}_1\underset k\otimes{\mathcal L}_2$ for some invertible sheaves ${\mathcal L}_1,{\mathcal L}_2$ on $G$ and $X$ respectively. We conclude that ${\mathcal L}$ is $G$-invariant.

(2) We shall need the following lemma:

\begin{lem} \label{lemaanterior} Let $G$ be a smooth and connected algebraic group acting on a variety $X$ and let $x\in X$ be a closed point such that $G$ acts faithfully on its orbit  $O_x$. Let   $D=\sum D_i$ be an effective Cartier divisor on $X$ whose irreducible components $D_i$ meet $O_x$ properly (i.e. $D_i\cap O_x\neq O_x$) and such that $x$ is an isolated point of $\cap D_i$. If ${\mathcal O}_X(D)$ is $G$-invariant, then $G$ is affine.
\end{lem}

\begin{pf} Let us denote $E=H^0(X,{\mathcal O}_X(D))$. It has no base points in $O_x$, because the base points subscheme is invariant  by $G$ and $D\cap O_x\neq O_x$. Let   $e\in E$ be such that the associated divisor is  $D$.  By Theorem \ref{SubInv0} and Definition \ref{definicionCoherente}, there exists a finite dimensional subspace  $V\subset E$  containing $e$ and stable under the action of  $G$.  Then we have a group  homomorphism  $\phi\colon G\to {\bf PGl}_k(V)$. The kernel of $\phi$ is finite: indeed, let  $H\subset G$ be a smooth connected subgroup  that acts trivially on ${\Bbb P} (V)$. Then $H$  acts trivially on  ${\Bbb P}(V^*)$; in particular,  $<e> $ (and then $D$) is fixed by $H$. Thus, all the irreducible components of  $D$ are fixed by $H$, so the reduced connected components of  $\cap D_i$ are fixed too. In particular,   $x$ is a fixed point under $H$. Since $G$ acts transitively on  $O_x$ and $H$ is normal in $G$, every point of  $O_x$ is fixed by $H$, i.e., $H$ acts trivially on  $O_x$. Since the action on $O_x$ is faithful, $H$ must be trivial. This proves that $\ker\phi$ is finite. Now,  ${\bf PGl}_k(V)$ is an affine group and  $\Img\phi$ is a closed subgroup (because  $\phi$ is a group homomorphism), so $\Img\phi$ is affine. Since $\phi\colon G\to \Img\phi$ is a finite morphism, $G$ is affine.
\end{pf}

 Since $G$ acts faithfully on $X$,  $\bigcap_{x\in G(k)}G_{x}=\{ 1\}$, where $G_{x}\subset G$ is the isotropy subgroup of $x$. Then  there exist rational no fixed   points $x_1,\dots ,x_m\in X$ such that $\cap_iG_{x_i}=\{ 1\}$. Moreover, since  $X$ is smooth, there exist divisors  $D_{ij}$ on $X$ satisfying: $D_{ij}$ contains $x_i$, meets  $O_{x_i}$ properly in $x_i$ and $x_i$ is an isolated point of  $\cap_jD_{ij}$ (both issues are a local problem and locally they are simple). Now, let us denote $ X^m=X\times \overset m\cdots\times X$,  $\pi_i\colon X^m\to X$   the i-th projection, $D=\sum_{ij}\pi_i^*D_{ij}$, ${\mathcal L}={\mathcal O}_{X^m}(D)$ and $x^{(m)}=(x_1,\dots ,x_m)$. One has: 

\item{(i)} ${\mathcal L}^n$ is $G$-invariant for some $n>0$, because  ${\mathcal O}_X(nD_{ij})$ is $G$-invariant for some $n>0$.

\item{(ii)} The irreducible components $\pi_i^*D_{ij}$ of $D$ meet  $O_{x^{(m)}}$ properly and $x^{(m)}$ is isolated in $\cap_{ij}\pi_i^*D_{ij}$.

\item{(iii)} $G$ acts faithfully on  $O_{x^{(m)}}$.  

One concludes by Lemma \ref{lemaanterior}. \end{pf}

\begin{cor}[\cite{Sumihiro2} Thm. 1.6, \cite{Knop} Thm. 2.4, \cite{Brion2} Thm. 2.14]\label{Linearizability}  Let $G$ be a connected and smooth affine algebraic group and let $X$ be a normal $G$-variety. For each invertible sheaf ${\mathcal L}$    on $X$   there exists $n>0$ such that ${\mathcal L}^n$ is  $G$-linearizable.\end{cor}

\begin{pf}  By Theorem \ref{grupoafin}, ${\mathcal L}$ is $G$-invariant. Moreover ${\mathcal L}^{\tau}$ is effective because $G$ is affine.  We conclude by Corollary \ref{torsion}.
\end{pf}

\begin{rem}\label{remarklinearizable} We have proved that Corollary \ref{Linearizability} is a consequence of Corollary \ref{torsion}. Proofs in \cite{Knop}, \cite{Sumihiro2} and \cite{Brion2}  use that $ {\rm Pic} (G) $ is finite, which is not true if $ G $ is not affine, and therefore they cannot be generalized to prove Corollary \ref{torsion}. That is, Corollary \ref{torsion} shows that the essential reason why Corollary \ref{Linearizability} holds is the invariance of the invertible sheaf and not the finiteness of the Picard group of $G$. Moreover, the finiteness of $ {\rm Pic} (G) $ is also a consequence of Corollary \ref{torsion}, as we will see now.\end{rem}

\begin{cor}[\cite{Grothendieck}, p. 5-21]\label{torsionful} If $G$ is an affine and smooth algebraic group (i.e. a linear group), then ${\rm Pic}(G)$ is finite.
\end{cor}

\begin{pf} We can suppose that $G$ is connected. By Proposition \ref{PicTrivialprop}, ${\rm Pic}(G)$ is a finitely generated conmutative group and then it is enough to prove that it is a torsion group. Let ${\mathcal L}$ be an invertible sheaf on $G$.  By Theorem \ref{grupoafin}, $\mathcal L$ is $G$-invariant and, by Corollary \ref{torsion}, ${\mathcal L}^n$ is fixed by $G$ (for some $n> 0$). One concludes by Remark \ref{miremark}.
\end{pf}

\medskip\section{Invariance of invertible sheaves by 1-dimensional deformations} 

\medskip

In this section we want to reduce the $G$-invariance of an invertible sheaf to the invariance by any smooth curve through the origin of $G$. More precisely:

\begin{thm}{\label{CurvoInv}} Let $G$ be a smooth and connected algebraic group acting on a variety $X$  and let ${\mathcal L}$ be an invertible sheaf on $X$. If ${\mathcal L}$ is invariant by the points of any locally closed smooth curve  $\C\hookrightarrow G$ passing through the origin, then there exists   $n> 0$ such that ${\mathcal L}^n$ is  $G$-invariant.\end{thm}

For the proof, we shall need some preliminary results.

\begin{lem} \label{torsionNucleoPic} Let $\pi\colon \overline X\to X$ be a finite and flat morphism of constant degree   $n$ and ${\mathcal L}$ an invertible sheaf on $X$ such that $\pi^*{\mathcal L}$ is trivial. Then  ${\mathcal L}^n$ is trivial.
\end{lem}

\begin{pf} One has
$${\bigwedge}^n_{{\mathcal O}_ X}(\pi_* {\mathcal O}_{\overline X})\simeq{\bigwedge}^n_{{\mathcal O}_{  X}}(\pi_*  \pi^* {\mathcal L})= {\mathcal L}^n \otimes_{{\mathcal O}_X} {\bigwedge}^n_{{\mathcal O}_{ X}}(\pi_*  {\mathcal O}_{ \overline X})  $$ hence ${\mathcal O}_X\simeq {\mathcal L}^n$.
\end{pf}

\begin{prop}\label{propoANterior} Let $X$ be scheme, $\tau\colon S\to {\bf Aut}_k(X)$ an automorphism and $\mathcal L$ an invertible sheaf on $X$. Let  $\pi\colon \overline S\to S$ be a finite and flat morphism of degree  $n$ and let us denote $\overline \tau=\tau\circ \pi\colon \overline S\to {\bf Aut}_k(X)$. Then

\begin{enumerate}
\item If ${\mathcal L}$ is fixed by  $\overline \tau$, then ${\mathcal L}^n$ is fixed by  $\tau$.

\item If $X$ and $S$ are varieties  and ${\mathcal L}$ is invariant by  $\overline \tau$, then ${\mathcal L}^n$ is invariant by   $\tau$.

\end{enumerate}
\end{prop}

\begin{pf} (1) If  ${\mathcal L}$ is fixed by  $\overline \tau$, then $(\overline\tau ^* {\mathcal L}_{\overline S}) \otimes_{{\mathcal O}_{X\times\overline S}}  ({\mathcal L}_{\overline S})^{-1}$ is  trivial. Since $$(\overline\tau ^* {\mathcal L}_{\overline S}) \otimes_{{\mathcal O}_{\overline S \times X}}  ({\mathcal L}_{\overline S})^{-1} = (\pi\times 1)^*\left[ ( \tau ^* {\mathcal L}_{  S}) \otimes_{{\mathcal O}_{ S \times X}}  ({\mathcal L}_{ S})^{-1}\right]$$ one has that $[( \tau ^* {\mathcal L}_{  S})  \otimes_{{\mathcal O}_{ S \times X}}  ({\mathcal L}_{ S})^{-1}]^n$ is trivial by Lemma \ref{torsionNucleoPic}. Hence ${\mathcal L}^n$ is fixed by  $\tau$.

(2) follows from (1) and Remark \ref{invislocfixed}.\end{pf}

\begin{prop}\label{propoANterior2} Let $G$ be a smooth and connected algebraic group acting on a variety  $X$  and let ${\mathcal L}$ be an invertible sheaf on $X$. Let $S$ be a variety and  $g\colon  S\to G$ an epimorphism. If  ${\mathcal L}$ is invariant by $ g$, then ${\mathcal L}^n$ is $G$-invariant for some $n> 0$.\end{prop}

\begin{pf}  Let $p$ be the generic point of  $G$. Since $g$ is surjective, the fiber of $p$ is non empty. By Hilbert's Nullstellensatz, there exists a point   $q\in  S$ in the fibre of $p$ whose residual field $k(q)$ is a finite extension of  $k(p)$. Let $S'$ be the closure of  $q$, $S'=\overline q$. It is clear that $g\colon S'\to G$ is generically finite, hence it is  finite over
an open $U$ of $G$. After localization, we may assume that  $g\colon S'\to U$ is finite and flat. By  Proposition \ref{propoANterior}, there exists  $n$ such that ${\mathcal L}^n$ is invariant by any point of  $U$. On the other hand,  ${\mathcal L}^n$ is invariant by any rational point of  $G$, since $g\colon  S\to G$ is surjective and then it is also surjective on rational points. By  Proposition \ref{AccionAccion0}, ${\mathcal L}^n$ is invariant by $g_0\cdot U$ for any   rational point $g_0\in G$. Since the open sets $g_0\cdot U$ cover $G$, we are done.\end{pf}

\begin{lem}\label{curvas} Let $G$ be a smooth and connected algebraic group. There exist smooth and connected locally closed curves  $C_1,\dots , C_n\hookrightarrow G$  passing through the origin such that $G=C_1\cdots C_n$.\end{lem}

\begin{pf} If $V\subset X$ is of the form $V=C_1\cdots C_n$, with $C_1, \dots ,C_n\hookrightarrow G$ smooth and connected locally closed curves passing through the origin $e\in G(k)$, then $V$ is an irreducible and constructible  set (i.e. a union of subvarieties of $G$); hence, the closure  $\overline V$ is an irreducible subvariety and   $V$ contains a dense open  subset $V_0$ of $\overline V$   (hence $dim\, V_0=dim\,  V=dim\,  \overline V$).

Let  $V$ be as above with maximal dimension. Let us see that $V$ contains a non empty open subset of  $G$.  In fact, if $dim\, V<dim\, G$, let us choose a point $x_0\in V_0$ and a smooth and connected curve $C\subset X$ passing through    $x_0$ and such that $C\cap \overline V$ is finite. One has $dim\,  C\cdot V=dim\, (x_0^{-1}\cdot C)\cdot V$ and, by maximality,   $dim\, (x_0^{-1}\cdot C)\cdot V=dim\, V$. Therefore $V_0$ is an open subset of $C\cdot V$, so they coincide locally at  $x_0$. This is not possible because $C\cap V_0$ is finite and $C\cdot V$ contains $C$.

To conclude, let us see that  $V\cdot V=G$. Let $x\in G(k)$ be any rational point; let us consider the inversion morphism  $i\colon G\to G$, $i(g)=g^{-1}$. One has that $x\cdot i(V)$ contains an open subset of $G$ (because $V$ does), and hence  $V\cap (x\cdot i(V))$ is not empty. Therefore, there exist rational points $v_1,v_2\in V(k)$ such that $v_1=x\cdot v_2^{-1}$, and then $x=v_1\cdot v_2\in V\cdot V$.\end{pf}

Now let us prove Theorem \ref{CurvoInv}.

\begin{pf} By Lemma \ref{curvas} there exist smooth curves  $\C_1,\dots ,\C_n\subset G$ passing through the origin such that $\C_1\cdots \C_n=G$, i.e., such that the multiplication morphism $g\colon \C_1\times\cdots\times \C_n\to G$ is surjective. By the hypothesis and  Proposition  \ref{AccionAccion0}, $\mathcal L$ is invariant by $g$; now, one concludes by Proposition \ref{propoANterior2}.\end{pf}

\section{Proof of  the Chevalley-Barsotti Theorem on algebraic groups}\medskip

The aim of this section is to give a proof of the Chevalley-Barsotti structure theorem on algebraic groups (Theorem \ref{Chevalley-Barsotti}).

\begin{prop}\label{extenMorfisPic} Let $X$ be a smooth variety,  $i\colon U\hookrightarrow X$ an open immersion. The morphism
$$i^*\colon {\mathbf Pic}^0(X)\to {\mathbf Pic}^0(U)$$
is an isomorphism over pointed smooth schemes; that is, for any pointed smooth scheme $(S,s_0)$, the map between pointed homomorphisms
$${\rm Hom} ((S,s_0), {\mathbf Pic}^0(X)) \to {\rm Hom} ((S,s_0), {\mathbf Pic}^0(U))$$
is bijective.\end{prop}

\begin{pf} By Lemma \ref{sucExac} one has that  $i^*$ is surjective and its kernel is  the pointed subfunctor ${\bf {Pic}}^0_{X-U}(X)$. Since $ {\mathbf {Pic}} _ {X-U} (X) $ is representable by a discrete scheme, we conclude that $ {\mathbf {Pic}} ^ 0_ {X-U} (X) = 0 $.\end{pf}

\begin{prop} \label{HomENAbel} Let $G$ be a smooth and connected algebraic group,  $\C$ a smooth curve. If $T\colon G\to {\mathbf Pic}^0(\C)$ is a pointed morphism, then it is a morphism of group functors.\end{prop}

\begin{pf}  By Proposition \ref{extenMorfisPic},  we may assume that  $\C$ is complete. Let  $\varphi \colon G\times G\to {\mathbf Pic}^0(\C)$ be defined by $\varphi (g,\overline g)=T(g\overline g)-T(g)-T(\overline g)$. One has to prove that   $\varphi $ is null. By hypothesis, it vanishes on  $\{ 1\}\times G$ and $G\times \{ 1\}$. Let $Z$ be the locus of points where   $\varphi $ vanishes. It is a closed subscheme, because  $\varphi $ is a morphism of schemes.

It suffices to prove that $\varphi$ vanishes on  $\C'\times G$ for any smooth and connected curve  $\C'\subset G$ passing through the origin; indeed, if this is the case and $Z\neq G\times G$, then $Z$ does not contain some smooth and connected curve  $\C''\subset G\times G$ meeting transversally $\{1\} \times G$ at $(1,1)$. Let $\pi_1\colon G\times G\to G$ be the first projection and $\C'=\pi_1(\C'')$. Then   $\C'$ is a smooth and connected curve (since it is aconstructible subset of the Zariski clousure $\overline{\pi_1(\C'')}$\,) 
through $\{ 1\}$  and $\C''\subset \C'\times G$. But $\C'\times G\subset Z$ by hypothesis and  we obtain the contradiction $\C''\subset Z$. 

Now, if $ \C '$ is complete, then $\varphi \colon \C'\times G\to {\mathbf {Pic}}^0(\C)$ is null because it is so over $ \C '\times \{1 \} $ and then, by the standard Rigidity Lemma, it is constant (null) on each fiber  $\C'\times \{ g\}$. 

If $\C'$ is not complete,   it suffices to prove that $\varphi \colon  \C'\times G\to {\mathbf {Pic}}^0\C$ extends to a morphism $\varphi\colon \overline {\C'}\times G\to {\mathbf Pic}^0(\C)$, where  $\overline {\C'}$ is the completion of $\C'$. That amounts to proving that the invertible sheaf $\mathcal L_{\varphi }$ on $\C'\times G\times \C$ corresponding to $\varphi$ extends to an invertible sheaf on $\overline {\C'}\times G\times \C$; equivalently, one has to prove that the morphism  $G\times \C\to {\mathbf Pic}^0(\C')$ defined by $\mathcal L_{\varphi }$ extends to a morphism $G\times \C\to {\mathbf Pic}^0(\overline {\C'})$. This follows from  Proposition \ref{extenMorfisPic}. 
\end{pf}

\begin{thm}{\bf (Chevalley-Barsotti)} \label{Chevalley-Barsotti} Let $G$ be a smooth and connected algebraic group over an algebraically closed field. There exists a unique smooth, connected, affine and normal subgroup  $H\subset G$ such that $G/H$ is an abelian variety. \end{thm}

\begin{pf} Let $H'$ be the subfunctor of the functor of points of $G$ consisting of those points that leave all the invertible sheaves of $G$ invariant (by translation), and let $H$ be the  reduced and connected component through the origin of $H'$ (i.e., $H={H'}^0$, see Section \ref{Notations}). It suffices to prove that $H$ is representable and that  $G/H$ is an abelian variety; indeed, $H$ is affine  by (2) of Theorem  \ref{grupoafin};  
moreover, if $N\subset G$ is a smooth, connected, affine and normal subgroup such that $G/N$ is an abelian variety then $N\subseteq H$ and $H/N$ is an affine, smooth and connected subgroup of the abelian variety $G/N$ and then $H/N=\{ e\}$, i.e.,  $N=H$.

Let us denote $L\colon G\times G\to G$ the left action, $L(g,g')=gg'$. For each invertible sheaf  $\mathcal L$ on $G$ and each smooth curve through the origin  $\C\subset G$, let  $$\aligned \varphi_{{\mathcal L},\C}\colon G&\to {\bf Pic}^0( \C)\\ g&\mapsto (L_g^*{\mathcal L}_S\otimes {\mathcal L}_S^{-1})_{\vert S\times \C}\endaligned$$  for any $S$-valued point $g\colon S\to G$. Let $\overline \C$ be the the completion of $\C$ and $\overline\varphi_{{\mathcal L},\C}\colon G\to {\bf Pic}^0(\overline \C)$ the unique extension of  $\varphi_{{\mathcal L},\C}$ such that $\overline\varphi_{{\mathcal L},\C}(1)=0$ (Proposition \ref{extenMorfisPic}). By Proposition \ref{HomENAbel}, $\overline\varphi_{{\mathcal L},\C}$ is a morphism of groups. 

Let $\overline H=\cap \ker \overline\varphi_{{\mathcal L},\C}$, which is representable. Let us see that $H$ is the reduced and connected component through the origin of $\overline H$. By definition, $H$ is a subfunctor of $\overline H^0$, so we have to prove that any invertible sheaf on $G$ is invariant by $\overline H^0$. By (1) of Theorem  \ref{grupoafin}, it is enough to prove that  $\overline H^0$ is affine. Let $D=\sum_iD_i$ be a divisor on  $G$ whose irreducible components $D_i$ pass through  $1$, meet  $\overline H^0$ properly at $1$ and their intersection contains $1$ as an isolated point. Let   ${\mathcal L}={\mathcal O}_G( D)$. For any  point $c$ of any smooth and connected curve $\C\subset G$ through $1$ and any point  $g$ in $\ker \varphi_{{\mathcal L},\C}$ one has ${\mathcal L}_{g\cdot c}\simeq {\mathcal L}_g\otimes {\mathcal L}_c$. In particular, ${\mathcal L}  $ is invariant by right translation by the points of any smooth curve of $\overline H^0$ passing through the origin. By Theorem \ref{CurvoInv}, ${\mathcal L}^n$ is right-$\overline H^0$-invariant (for some $n> 0$). By Lemma \ref{lemaanterior},  $\overline H^0$ is affine.

To conclude, let us see that $G/H$ is an abelian variety. There exist a finite number of curves  $\C_1,\dots ,\C_n$ such that $\overline H=\cap_{i=1}^n \ker \overline\varphi_{{\mathcal L},\C_i}$, so  $G/\overline H\hookrightarrow  \prod_{i=1}^n {\bf Pic}^0(\overline \C_i)$ is an abelian variety. Since $G/H\to  G/\overline H$ is finite, we conclude. 
\newcommand{\potato}{{\it Unicity:} Let  $H_2$ be another smooth and connected affine subgroup. The image of $H_2$ by $G\to G/H$ is an abelian subvariety of $G/H$ and it is affine (because it is a quotient of  $H_2$). Hence this image is trivial, and then   $H_2\subset H$. Analogously, $H\subset H_2$.}\end{pf}

\begin{rem} The proof of this theorem gives a construction of the subgroup $G_{\rm aff}\subset G$ as the reduced and connected component through the origin of the subscheme of points of $G$   leaving invariant (under traslation) any invertible sheaf on $G$ (an ample invertible sheaf would suffice). 
\end{rem}

\end{document}